# An Approach to Development: Turning Education from a Service Duty to a Productive Tool

Pooya Alinian[1], Raziyeh Mohammadi[2], Azadeh Parvaneh[3], Ali Rejali[4]

[1] Department of Mathematical Sciences, Isfahan University of Technology, Isfahan, Iran.
Email: pooyaalinian76@gmail.com
ORCID: https://orcid.org/0000-0003-3347-7406

[2] Department of Mathematical Sciences, Isfahan University of Technology, Isfahan, Iran.
Email: razieh.mohammadi@math.iut.ac.ir
ORCID: https://orcid.org/0000-0002-5226-1734

[3] Department of Statistics, University of Isfahan, Isfahan, Iran.
Email: azadee.parvanee@yahoo.com
ORCID: https://orcid.org/0000-0002-0594-3601

[4] Department of Mathematical Sciences, Isfahan University of Technology, Isfahan, Iran.
Email: a rejali@cc.iut.ac.ir
ORCID: https://orcid.org/0000-0002-5747-6646
*(Corresponding author)*

**Abstract**

Recent economic developments of countries like Japan, Korea, and Singapore, as a result of improvement in the quality of their education, show that having a high-quality education may lead to economic growth. In this article, using some statistical methods, we argue that high quality education can change the economy towards higher growth. Therefore, for the development of the country, one should think about how to improve its education. One of the effective ways to improve the quality of education is to increase the efficiency of teachers and attract talented people to teaching positions. Research shows that raising teachers' salaries, along with a proper quality improvement program, can help facilitate this process.

# Introduction

Improving the quality of the workforce is one of the essential strategies for achieving economic growth and ultimately sustainable development. The quality of the workforce is guaranteed by investing in human education. In particular, well-trained future generations will guarantee the quality of the workforce. The experiences of developed countries show that explaining economic growth rates, only based on increasing physical investments and working population, is inadequate and inaccurate, and another factor has contributed to the economic growth of these societies, as well. This factor, known as the "residual factor", is the main cause for the increase in the productivity of capital and manpower. Economists believe that the residual factor, which explains an important part of the economic growth in developed countries, is directly or indirectly related to the improvement of their education (Emadzadeh, 1992).

In less developed countries, including Iran, human resources despite their tremendous impact on sustainable development, have been ignored. The education officials look at education as a public service, and the executive approach to it is similar to providing health care, transportation or utilities, while it is necessary for the decision-makers and also all the citizens to come to a common understanding that any expenditure on education in a country will lead to added values, and removing the barriers for the qualified education makes the country able to achieve sustainable development. In other words, achieving sustainable development requires changing the look in education from the place of public service to the most important productive tool, which is producing the greatest capital needed for human growth (Rejali, 2009).

Since the mechanism for turning education from a service duty to a productive tool requires careful scientific research and investing on this subject is a long-term activity, countries' political decision-makers do not easily accept to spend on this procedure unless they become sure of its benefits. To achieve this goal, we intend to investigate the role of education for sustainable development. Accordingly, we show that improving the quality of education leads to economic growth. Some reasons for the success of the education systems of the important countries like Japan, Korea, and Singapore, which are countries that have achieved sustainable development by improving the quality of their education, are expressed as a model for less developed countries, including Iran. Then, we mention Iran's disappointing condition in education as a developing country. Finally, based on the results of similar studies, the most effective factor in improving the educational structure and learning process of students is introduced, and one of the most important strategies for achieving it, will be presented.

# Main Problem: Is Education a Service Duty of the Government and the Parents, or Should It Be Considered as a Productive Tool?

Looking at education as a public service and considering it as a government duty in the constitution of many countries is a short-sighted view that now dominates the education of many countries, unfortunately, and nowadays this service, with the help of so-called nonprofit organizations, various foundations, charities, and even the families of the children, is delivered to these societies like a luxury commodity and people must be debtors to these institutions for these kinds of services. But education is not a service work, it is an important productive tool, and actually one of the most important pillars of countries' productions. Countries with abundant economic resources but lacking the intellectual and creative manpower are always colony and have to share the revenues of their products with others so that they can develop them if they like (Rejali, 2009). In this paper, in order to show that education is a productive tool, we first show that there exists a significantly positive statistical relationship between improving education quality and higher economic growth.

# The Relationship between Education and Economic Growth

*Literature Review*

Since the late 1980s, most of the macroeconomic attention has focused on determining factors for a long-term economic growth including education (Barro, 2001). The empirical framework summarized in (Barro, 1997) shows that the "growth rate" is a function of "per capita product" and "long-run level of per capita product". When human capital is considered, per capita product would be generalized to include the levels of physical and human capital. In some theories, the growth rate increases with the ratio of human to physical capital (Barro, 2001). "Human capital is productive wealth embodied in labor, skills and knowledge."[1]

Barrow (Barro, 2001) emphasizes on the role of education in long-term economic growth and shows that:

- "Growth is positively related to the starting level of average years of school attainment of adult males at the secondary and higher levels."

- "Growth is insignificantly related to years of school attainment of females at the secondary and higher levels. This result suggests that highly educated women are not well utilized in the labor markets of many countries."

---

[1] See: https://stats.oecd.org/glossary/detail.asp?ID=1264

- "Growth is insignificantly related to male schooling at the primary level. However, this schooling is a prerequisite for secondary schooling and would, therefore, affect growth through this channel."

- "Education of women at the primary level stimulates growth indirectly by inducing a lower fertility rate."

*Methodology and Results of a New Trial for Estimating the Relationship between Education and Economic Growth*

As another study on the relationship between education and economic growth, we examined the impact of quality of education on economic growth in 137 countries whose information were available. The necessary information for checking these issues were gathered from the Global Competitiveness Report 2017-2018[2]. The primitive benchmarks that have been considered as the default for measuring the quality of education include:

1. Quality of primary education;
2. Primary education enrollment rate;
3. Secondary education enrollment rate;
4. Tertiary education enrollment rate;
5. Quality of the education system;
6. Quality of math and science education;
7. Quality of management schools;
8. Internet access in schools;
9. Local availability of specialized training services;
10. Extent of staff training.

According to a factor analysis, these 10 variables are divided into 4 following factors (see Figure 1):

$V_1$: Quality of education (Quality of primary education, Quality of math and science education, Quality of management schools, Internet access in schools, Local availability of specialized training services, Extent of staff training);

---

[2] The Report and an interactive data platform are available at: www.weforum.org/reports/theglobal-competitiveness-report-2017-2018

$V_2$ : Primary education enrollment rate;

$V_3$ : Enrollment rate in high school and higher education (Secondary education enrollment rate, Tertiary education enrollment rate);

$V_4$ : Quality of the education system.

We then consider the following criterion for measuring the quality of education, which is a weighted average of the four factors obtained by the factor analysis:

$$\text{QA} = \frac{0.59}{0.87}V_1 + \frac{0.12}{0.87}V_2 + \frac{0.09}{0.87}V_3 + \frac{0.07}{0.87}V_4.$$

Now the relationship between education quality and economic growth is examined by studying the statistical relationship between QA and GDP per capita in US$ obtained from the Global Competitiveness Report 2017-2018. To end this part of the study, we examined the effect of the increase in QA (as a predictor variable of the model) on GDP per capita in US$ (as the response variable of the model), by fitting an appropriate regression model.

To increase the accuracy of the procedure, we added the following predictor variables in the model:

$V_5$ : Inflation, GDP deflator (annual %);

$V_6$ : International migrant stock, total;

$V_7$ : Terms of trade adjustment (constant LCU);

$V_8$ : Trade (% of GDP);

$V_9$ : Population, total;

$V_{10}$ : Export volume index (2000 = 100);

$V_{11}$ : Total natural resources rents (% of GDP).

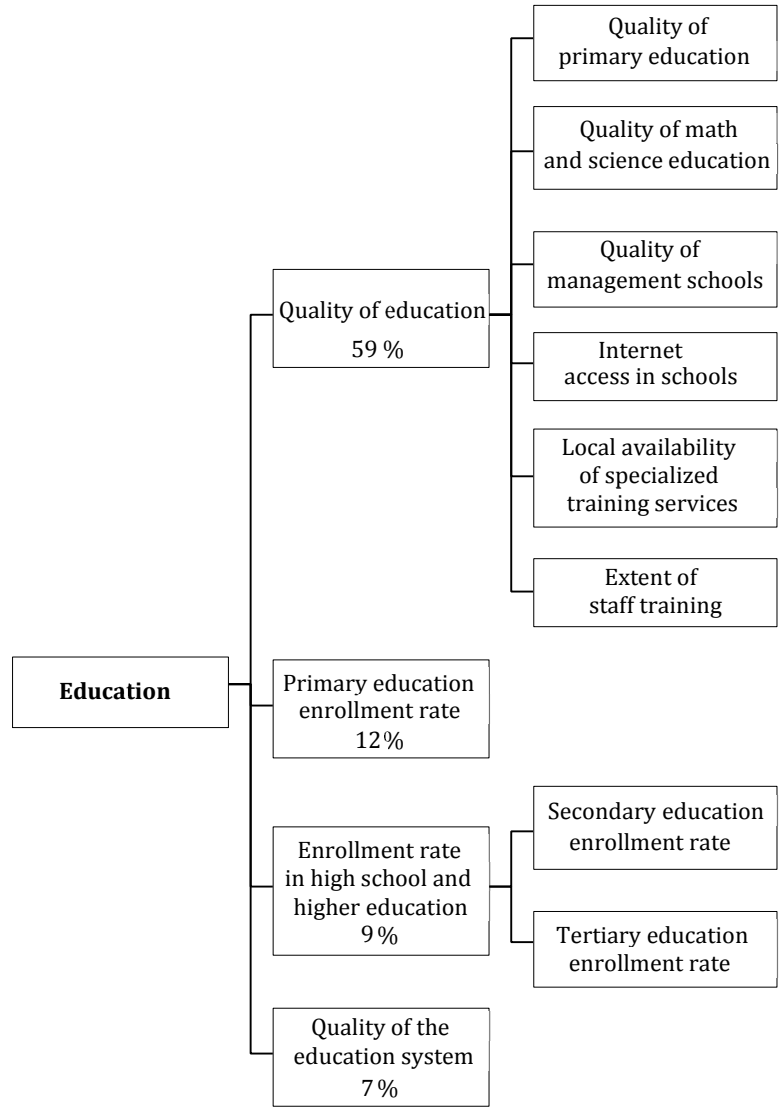

Figure 1: *The result of the factor analysis. Each number next to the agent represents the percentage of information that this factor justifies the related primary variables.*

These random variables have been introduced in the World Bank website[3], so for more information, the reader can visit its website. We also collected the observations related to the above variables from the World Bank website, all of which, except $V_6$, are for the year 2018, and only for $V_6$, due to the lack of data in 2018, 2015 data was used. Moreover due the fact that the World Bank information does not include any information for Taiwan and China, separately, while the World Economic Forum provides distinct information for the QA and GDP per capita in US$ of them, so we considered China's observations for these two variables as a weighted average of Taiwan and China data with weights related to both countries' population, which were reported by the Global Competitiveness Report 2017-2018.

---

[3] https://www.worldbank.org

But we estimated the missing data of other countries using the following two approaches:

1. If some data for a country were missed in one year, while there was some suitable information from previous years, the missing data has been forecasted by using an ARIMA time series model (Ansley & Kohn, 1983). The missing values which were estimated by this method are observations for $V_8$ of Iran, Laos, Tajikistan, Tanzania and Venezuela, as well as the observation of $V_5$ for Iran. We could not use this method for forecasting the rest of the missing values because of lacking suitable information for previous years or some increase in forecasting variance.

2. To estimate the rest of the missing values, the random forest method is used. In this method, several trees have been formed by using the available variables, and missing values were being predicted based on these trees (Tang & Ishwaran, 2017). The missing data which were predicted by using this method are the observations for $V_5$ of Yemen, $V_7$ of Andorra, Azerbaijan, China, Iran, Laos, Qatar, Seychelles, Tajikistan, Tanzania, Trinidad and Tobago, Tunisia, Venezuela, Yemen and Zambia, $V_8$ of Andorra, Trinidad and Tobago and Yemen, $V_{10}$ of Montenegro and Switzerland, and $V_{11}$ of Iran and Venezuela.

Finally, in order to be able to compare that which one of the variables has the most effect on the response variable, we standardized all the predictor variables, and then consider a linear model with the response variable

$$V = \sqrt{\text{GDP per capita in US\$}}$$

against predictor variables SQA, $SV_5$, $SV_6$, $SV_7$, $SV_8$, $SV_9$, $SV_{10}$ and $SV_{11}$, where SQA and $SV_i$, $5 \leq i \leq 11$, are the standardized QA and the standardized $V_i$, respectively. The result of fitting this model has been mentioned in Table 1. According to the table, SQA has the most effect on the response variable in comparison with the other significant variables, and its effect is positive. Also, based on the adjusted R-squared, the percentage of justification of the response variable by this model is an acceptable and valid value. Generally, this result shows that, in addition to the positive effect of education quality on economic growth, the quality of education has the most effect on economic growth compared to other variables studied in the model which cover different aspects from economic to demographic values.

Table 1: *The output of the regression with response variable V and explanatory variables SQA, $SV_5$, $SV_6$, $SV_7$, $SV_8$, $SV_9$, $SV_{10}$ and $SV_{11}$*

|  | **Coefficient** | | |
| --- | --- | --- | --- |
| **Variable** | Estimate | Std. Error | P-Value |
| Intercept | 101.03 | 3.10 | 0.00 |
| SQA | 54.89 | 3.64 | 0.00 |
| $SV_5$ | 5.04 | 3.34 | 0.13 |
| $SV_6$ | 7.21 | 3.40 | 0.04 |
| $SV_7$ | -1.33 | 3.24 | 0.68 |
| $SV_8$ | 12.12 | 3.48 | 0.00 |
| $SV_9$ | -6.72 | 3.29 | 0.04 |
| $SV_{10}$ | -0.78 | 3.27 | 0.81 |
| $SV_{11}$ | 0.78 | 3.34 | 0.81 |

denotes significance coefficients at 5%
Multiple R-squared: 0.75   and   Adjusted R-squared: 0.74

## Japan, Korea, and Singapore's Development after Improving Their Educational Situation

Japan, Korea, and Singapore are the three countries with few natural resources, but they have experienced rapid and generally sustained their economic growth, since 1960s. In recent years, they have increasingly sought a policy of combining education and human resources planning with economic policies for sustainable development. Trends in International Mathematics and Science Study (TIMSS)[4] is a type of study that has monitored trends in mathematics and science achievement of participated countries, every four years and at the fourth and eighth grades. The result of the Pearson correlation coefficient between TIMSS average mathematics scores of 8th-grade students in 2015 and QA with a p-value of 0.00 and estimated correlation of 0.60 shows a relatively high (positive) significant correlation. This means that TIMSS scores can be an alternative variable for QA in investigating the quality of education which has a positive impact on economic growth. Japan, Korea, and Singapore have always been among the best in TIMSS studies. One possible explanation for these achievements could be their attitude toward their human resources, because they put a strong emphasis on their education (Yang & Yorozu, 2015).

"In Singapore, for example, the government believes that the only way the country can compete with other developed countries is through heavy investment in education and the continuous upgrading of labor skills through training and development programs" (Yang & Yorozu, 2015, page 13). "The two largest budget items of the government expenditure are defense and education. Education for primary, secondary, and tertiary levels is mostly supported by the state" (Kaur, 2014, page 1).

---

[4] https://timssandpirls.bc.edu

In the education system of Singapore, which has always succeeded in international tests such as TIMSS (Table 2 shows that Singapore has been ranked first in TIMSS for most years), Singapore's teachers are often appreciated as the reason for Singapore's successes at international tests so that Asia Society (AS) believes that Singapore has one of the highest quality education in the world, due to the quality of its teachers. Singapore's policy for its teachers is that they carefully select, educate, develop, and properly reward. According to the Grattan Institute in Australia, a strong culture of teacher education, collaboration, counseling, and professional development, etc., explains the superiority of Singapore's education (King, 2016). A report by the American Aspen Institute describes Singapore as a model for teacher development (Scalfini & Lim, 2008). In fact, the Singapore Ministry of Education hires prospective teachers from one-third of the best high school graduates (King, 2016). Moreover, Lee Kuan Yew, the founding father of modern Singapore, talks about teacher influence:

> *Outside the influence of their parents and their homes, the most important influence is the teacher and the school. And so it is no exaggeration to say that our 10,600 teachers in all our schools constitute the most influential group of 10,600 people anywhere in Singapore (Lee, 1959).*

"Lee recognised that emphasis must be placed on making teaching an attractive profession if we were to attract and recruit enough good people for the profession. He made it clear that a good pay is a necessary, albeit insufficient means, in achieving the end" (Tan et al., 2017, page 33). Also, Tony Tan Keng Yam, a Singaporean politician, says:

> *The basic solution to the problem of attracting high-quality people to join the teaching profession is to pay them sufficiently (Collins et al., 1991).*

Japan has also made remarkable achievements in recent decades, that many researchers search the secret of it in Japanese education culture. Japanese officials are particularly paying attention to elementary education, so that the government of Japan spends most of its GDP on preschool and elementary education (Bakhtiari, 2008). Arinori Mori, an architect of Japan's modernization, wrote in 1885: "Our country must move from its third-class position to second class, and from second class to first, and ultimately to the leading position among all countries of the world. The best way to do this is [by laying] the foundations of elementary education" (Lewis, 1995, page 10)

Another successful country in international tests such as TIMSS is Korea. From July 15th to July 17th, 2012, the United States National Commission on Mathematics Instruction and Seoul National University held a joint workshop on mathematics, teaching, and curriculum. The workshop was supposed to address matters related to mathematics teaching and curriculum in each country. Among the features of successful education system discussed in the workshop was its attention towards the teachers. In Korea, there are two categories of "first-class teachers" and "second class teachers" whose requirements for becoming a first-class teacher are doing extra study and passing an exam. Korean teachers earn points or credits that go toward their promotions and the chance to catch higher positions during their careers (National Academies of Sciences, Engineering, and Medicine, 2015).

As a summary of what was said, it is worth mentioning that all three education systems mentioned above, have been focused on their teachers and increasing their teachers' efficiency for improving their quality of education.

## The Situation of Education in Iran

In previous sections, it has been shown that improving the quality of education is effective in achieving economic growth which is an important aspect for sustainable development. But how is the quality of education in Iran? Very low scores in TIMSS (see Table 2) shows the existence of a crisis in Iranian schools' education situation.

There are many different possible reasons for the low quality of education in Iran. Some of the most important of them include (Parvaneh & Rejali, 2017):

- Low position for the teaching job;

- Underestimating education by the decision-makers;

- Expanding the domain of preparation process for university entrance exams on one hand to the elementary schools (for preparing the students for entrance examinations of specified schools such as the so-called gifted students' schools) and on the other hand to the graduate students at doctoral level;

- Unemployment and a sharp increase in the number of students;

- Teacher training issues, establishing of Farhangian University (which is responsible for training of all the teachers) without suitable preparation, and reduction of teachers' position, as well as changing school curriculum rapidly, without preparing the teachers to teach new materials.

Table 2: *Iran's TIMSS Results in Mathematics in comparison with other countries*

| Result | Year | | | | | |
|---|---|---|---|---|---|---|
|  | 1995 | 1999 | 2003 | 2007 | 2011 | 2015 |
| **4th Grade** | | | | | | |
| Iran's rank (number of participants) | 25 (26) | * | 22 (25) | 28 (36) | 43 (50) | 42 (48) |
| Iran's average score | 429 | * | 389 | 402 | 431 | 431 |
| International average | 529 | * | 495 | 473 | 491 | 507 |
| Top rank country | Singapore | * | Singapore | Hong Kong | Singapore | Singapore |
| Top average score | 625 | * | 594 | 607 | 606 | 618 |
| **8th Grade** | | | | | | |
| Iran's rank (number of participants) | 38 (41) | 33 (38) | 33 (45) | 34 (48) | 32 (42) | 29 (37) |
| Iran's average score | 428 | 422 | 411 | 403 | 415 | 436 |
| International average | 513 | 487 | 466 | 452 | 467 | 487 |
| Top rank country | Singapore | Singapore | Singapore | Taiwan | South Korea | Singapore |
| Top average score | 643 | 604 | 605 | 598 | 613 | 621 |

* In 1999, TIMSS was not held in 4th grade.

Source: www.nces.ed.gov/TIMSS

Ignoring teachers and their minimum needs, as well as their low values in society, reduced the value of science and teachers among students and even society. Most teachers have to choose a second job due to the need for more income, so they cannot have enough study in their spare time to be ready for their classes and students (Parvaneh & Rejali, 2017), and this is happening, while improving the quality of education can lead to economic growth in addition to cultural growth. Almost all less developed countries face similar problems in their education system and need to pay more attention to education and teachers.

## Increasing the Quality of Education by Increasing Teachers' Efficiency

*Literature Review*

In literature, different ideas for enhancing the quality of education have been investigated. For example, there have been many studies on the idea of reducing classroom size. The researchers' findings indicate that the policy of reducing the number of students in a classroom, although it is generally positive for both teachers and students, either has a negligible effect on students' performances or it is ineffective. In a meta-analysis of nearly 60 studies, Hanushek found that less than 15% reported a significant positive effect on the policy of reducing the number of students in a classroom on academic performances (Hanushek, 1997, 2003). Challenging the aforementioned meta-analysis method, Krueger (1999) presented a modified version and showed that only one-third of the studies in this modified meta-analysis have shown a significant positive effect on the policy of reducing the number of students in a classroom (Krueger, 1999,2002, 2003). Thus, most studies (between 66% and 85%) which have reviewed by Hanushek and Kruger analysis show that decreasing the number of students in a classroom does not affect students learning (Jensen, 2010a). Therefore, given the high cost of making classes with small sizes, countries should pursue a more appropriate policy to enhance the quality of their education.

*Methodology and Examination of a New Trial for Estimating the Impact of a Factor on the Quality of Education*

Another policy that many believe that may have a positive impact on education is the expenditure on education. Findings of a research (Parvaneh et al., 2019) indicate that the policy of increasing the education budget is also ineffective in many situations. In order to investigate this claim again, we tried to fit a random intercept model for "mathematics TIMSS scores for 4th grade"[5] ($Y$) as the response variable and "government expenditure per student, primary (% of GDP per capita)"[6] as the predictor variable, to account the education budget's effects. However, in order not to lose information, next to the TIMSS score of each year, for the expenditure, we used the average available information between the last four years. First, we considered the normal distribution for both conditional errors and random effects, but the diagnostic measures showed that the normal distribution is not appropriate for the conditional distribution of *Y* given the random effect. Figure 2 shows the plot of the conditional residuals of this model. Based on the figure, it seems that t and Laplace distribution may provide a better fitting, due to their fatter tails (Laplace distribution is sharper at the peak as well). So, we considered normal, t and Laplace distributions for the conditional errors and tried to fit them through the Bayesian analysis available in WinBUGS software.

---

[5] The data set is available at: www.nces.ed.gov/TIMSS (also at: https://timssandpirls.bc.edu)
[6] The data set is available at: https://data.worldbank.org

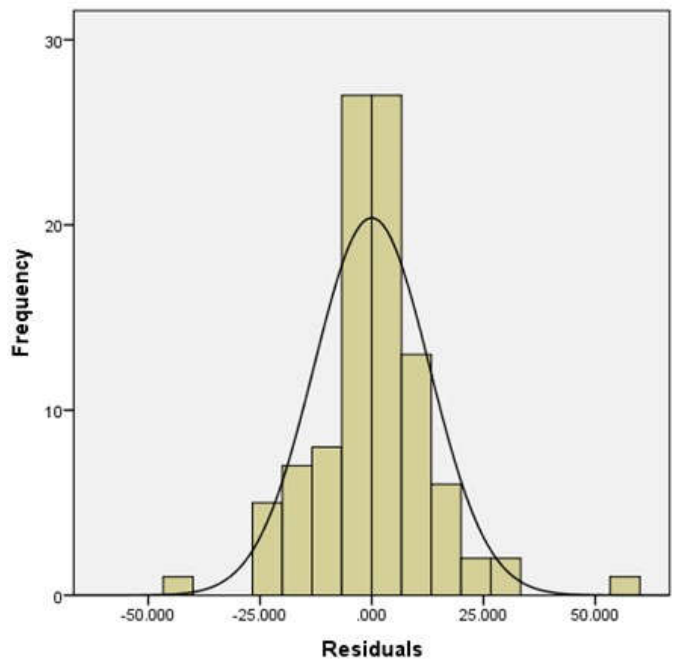

Figure 2: *Histogram of conditional residuals of the normal model*

For applying the Bayesian analysis, we considered noninformative normal distribution $N(0,10^6)$ for fixed effects $\beta_0$ and $\beta_1$, inverse Gaussian $IG(0.1,0.1)$ for the scaled parameters (this distribution is a noninformative prior with mean 1 (Sahu et al., 2003), and exponential distribution $Exp(0.1)$ truncated at 2 for degrees of freedom (this truncated guarantees the finiteness of the variance of the associated t distribution (Sahu et al., 2003). Also, we considered 10000 iterations after discarding the first 5000 iterations as burn-in to make the inference. We also checked the convergence through autoregressive, kernel density and other plots of parameters. Table 3 shows that the resulting estimated parameters including mean, standard deviation and the 95% credible intervals of parameters. Moreover, we report the deviance information criteria (DIC) that is a popular measure for comparing models in Bayesian analysis (Spiegelhalter et al., 2002). We note that a model with lower DIC is better fitted. Based on the DIC values of the three models, we can conclude that the model with Laplace distribution for conditional errors is more efficient than other models. However, in all three models, the confidence interval of $\beta_1$ contains 0, and accordingly, we can conclude that the expenditure on education does not have any effect on TIMSS scores or more precisely on a better quality of education.

Table 3: *Posterior means, standard deviation (SD) and 95% probability intervals for the parameters under three random intercept models for examining the effect of expenditure on education upon TIMSS scores*

| Model | DIC | Parameter | Mean | SD | 2.5% | 97.5% |
|---|---|---|---|---|---|---|
| Normal | 888.1 | $\beta_0$ | 485.1 | 12.52 | 462.2 | 513.8 |
| | | $\beta_1$ | 1.374 | 0.813 | -0.1879 | 3.017 |
| Laplace | 786.2 | $\beta_0$ | 477.4 | 19.85 | 434.1 | 509.3 |
| | | $\beta_1$ | 0.8845 | 0.6674 | -0.3901 | 2.25 |
| t | 838.2 | $\beta_0$ | 485.5 | 11.33 | 463.0 | 506.8 |
| | | $\beta_1$ | 0.9071 | 0.6386 | -0.3707 | 2.19 |

*Results and Discussion*

Undoubtedly, spending on education is necessary for increasing the quality of education, and our study only shows that raising the education budget, without proper planning on how it should be distributed, will not be so effective. More precisely, this policy requires further exploration that increasing the amount of spending on what part of education will be more effective. But what is an effective policy from the perspective of studies? Teachers' efficiency on the performance of students is more effective than any other plan or policy in education (Aaronson et al., 2007; Hanushek et al., 1998; Hanushek et al., 2005; Leigh & Ryan, 2010; Rockoff, 2004). There is enough evidence that the quality of a teacher influences the learning process. Quoted from (Jensen, 2010b), some studies in this area include:

- Hanushek (1992) estimated that the difference in results between students with a weak teacher and those with a well-prepared teacher could be as high as one year.

- Jordan et al. (1997) found that students who had three effective teachers successively, after three years, were 49% higher in school evaluations than students who had ineffective teachers.

- Sanders and Rivers (1996) found that students who had been training by effective mathematics teachers for three consecutive years had achieved approximately 50% higher scores than students who had similar mathematics scores at first but then had been training by less efficient teachers for three consecutive years.

- While having an effective teacher can help to achieve more than expected student results, an effective teacher may not fully make up the deficiency that the students have had in training by less effective teachers. Sanders and Rivers (1996) found that the students who, after being trained by some ineffective teachers, were trained by a highly effective teacher, although made more progress than anticipated, they have still not been compensated for their ineffective teachers.

Studies report some positive effects of increasing teachers' salary on their efficiency:

- A study by Parvaneh et al. (2019) shows that there is a significant positive correlation between teachers' salaries and TIMSS scores.

- A study by Loeb and Page (2000) shows that a 50% increase in teachers' salaries reduces high school dropout rates by more than 15% and increases college enrollment rates by about 8%. They also explain why previous studies have failed to produce systematic evidence that teacher wages affect student outcomes: "The identification of teacher wage effects has been based on cross sectional variation that reflects both supply and demand factors. Because alternative labor market opportunities and other school district characteristics vary across districts, the supply of teachers in quality-wage space also varies across districts" (Loeb & Page, 2000, page 406).

- A study by Figlio (1997) reveals that a significant relationship between teacher salary and teacher quality within local teacher labor markets exists. In this study, by considering "graduation from the highly selective undergraduate institution" and "majored in math, science, engineering or computer science" for the teacher qualifications, he has proved that the more teachers pay attracts the more qualified teachers in the metropolitan area of United States, and even in a particular metropolitan area, those school districts which pay higher salaries manage to employ better-qualified teachers.

Finally, we note that the policy of increasing the teachers' salary can significantly make efficiency and attract capable manpower if a detailed plan is arranged on how to implement it.

## Conclusion

The following aspects were considered in this article, namely

*Step 1.* Study of the impact of improving the quality of education on economic growth: Our results show that increasing the quality of education is highly effective.

*Step 2.* The implementation of three successful education systems that have achieved sustainable development through careful attention to their education and teachers shows that attention to the quality of education and the quality of teachers has a strong impact on sustainable development.

*Step 3.* The state of education in Iran is considered as an example of a developing country: Observations show that the status of education in Iran is low, and politicians should think about ways to improve the quality of education in the country.

*Step 4.* The research conducted in the article shows that the most effective way for increasing the quality of education is to increase the efficiency of teachers. And increasing their wages is one of the most effective ways for achieving this goal.

In particular, Figure 3 shows the relationship achieved in this study.

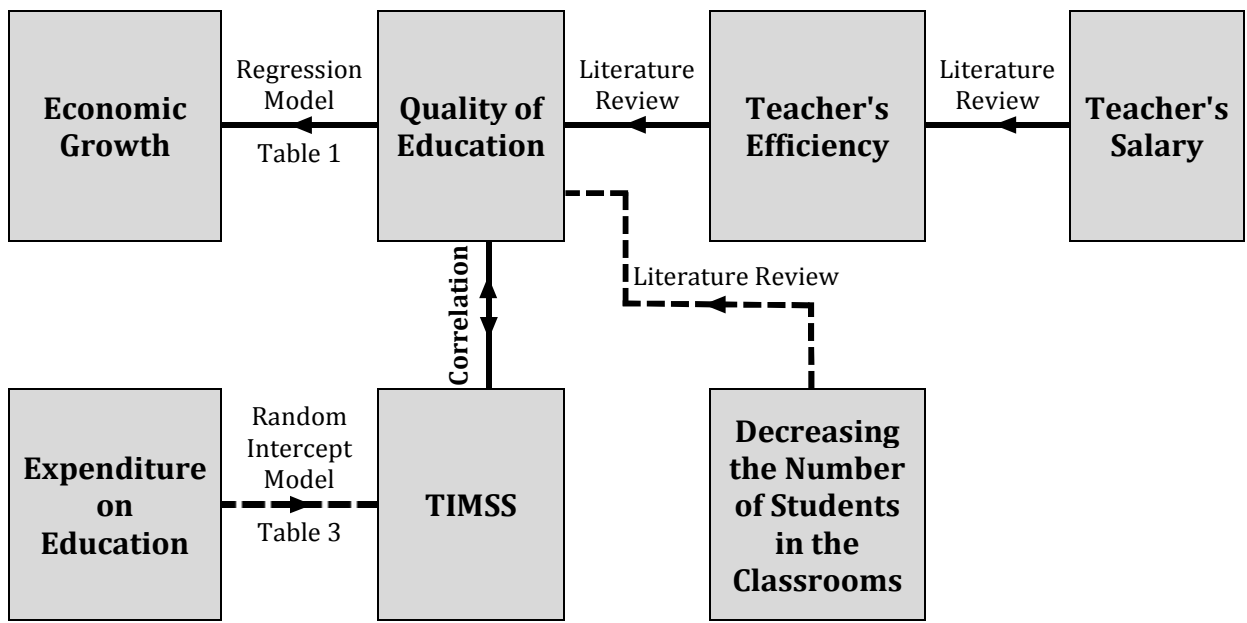

Figure 3: *The relationship between the various factors examined in this paper. The oriented line shows the existence of a significant direct effect, and the oriented dashed line indicates a lack of that.*